\documentclass{ifacconf}
\usepackage{graphicx}      
\usepackage{natbib}        

\newtheorem{theorem}{Theorem}[section]
\newtheorem{lemma}[theorem]{Lemma}

\newtheorem{result}[theorem]{Result}
\newtheorem{definition}{Definition}

\newtheorem{ass}{Assumption}

\usepackage{amsmath, amssymb}
\usepackage{amsfonts}
\usepackage{graphicx}
\usepackage{enumerate}
\usepackage[font = small]{caption}
\usepackage{ifthen}
\usepackage{multicol}
\usepackage{float}
\usepackage[usenames, dvipsnames]{color}
\usepackage[lofdepth,lotdepth]{subfig}
\graphicspath{{figures/}}
\usepackage{url}

\newcommand{\Ll}{L} 
\newcommand{\Gg}{\mathcal{G}} 
\newcommand{\Vv}{\mathcal{V}} 
\newcommand{\Ee}{\mathcal{E}} 
\newcommand{\Nn}{\mathcal{N}} 
\newcommand{\bll}{\bar{\lambda}_1} 
\newcommand{\nth}{$n^{\text{th}}$ }
\newcommand{\ddt}{\frac{\mathrm{d}}{\mathrm{d}t}}

\newcommand{\hinf}{$\mathcal{H}_\infty$ }

\newcommand{\grev}{grounded eigenvalue} 

\usepackage{array}
\newcolumntype{L}[1]{>{\raggedright\let\newline\\\arraybackslash\hspace{0pt}}m{#1}}
\newcolumntype{C}[1]{>{\centering\let\newline\\\arraybackslash\hspace{0pt}}m{#1}}
\newcolumntype{R}[1]{>{\raggedleft\let\newline\\\arraybackslash\hspace{0pt}}m{#1}}
\usepackage{subfig}

 \usepackage{tikz}
 \usetikzlibrary{decorations.pathmorphing,shadings}
 \usetikzlibrary{plotmarks}
 \usepackage{pgfplots}
\usepackage{tkz-graph}

\definecolor{gray3}{rgb}{0.75, 0.75, 0.75}
\definecolor{gray2}{rgb}{0.5, 0.5, 0.5}
\definecolor{gray1}{rgb}{0.25, 0.25, 0.25}
\definecolor{gray0}{rgb}{0.15, 0.15, 0.15}

\pgfplotscreateplotcyclelist{mygraylist}{%
black, densely dashed\\%
gray1\\%
gray3\\%
gray2\\%
gray2,dashed\\%
gray3,dashed\\%
}

\pgfplotscreateplotcyclelist{myothergraylist}{%
black, densely dashed\\%
gray1\\%
gray2\\%
gray3\\%
gray1, dashed\\%
gray3,dashed\\%
}

\definecolor{emmagreen1}{rgb}{0, 0.5, 0.1}
\definecolor{emmaorange1}{rgb}{0.99, 0.5, 0}
\definecolor{emmablue1}{rgb}{0, 0.25, 0.5}
\definecolor{emmaskyblue1}{rgb}{0.4, 0.8, 0.95}
\definecolor{emmapurple1}{rgb}{0.5, 0.25, 0.6}

\pgfplotscreateplotcyclelist{mycolorlist3}{%
emmagreen1\\%
emmaorange1\\%
emmablue1\\%
emmapurple1\\%
emmaskyblue1\\%
}

\usetikzlibrary{arrows}

\newboolean{showcomments}
\setboolean{showcomments}{true}

\newcommand{\rick}[1]{\ifthenelse{\boolean{showcomments}}
{\textcolor{blue}{Rick says: #1}}{}}
\newcommand{\emma}[1]{\ifthenelse{\boolean{showcomments}}
{\textcolor{red}{Emma says: #1}}{}}
\newcommand{\maria}[1]{\ifthenelse{\boolean{showcomments}}
{\textcolor{cyan}{Maria says: #1}}{}}
\newcommand{\todo}[1]{\ifthenelse{\boolean{showcomments}}
{\textcolor{Green}{To do: #1}}{}}

\theoremstyle{nonumberplain}%
\usepackage[font=footnotesize]{caption}

\begin{document}
\begin{frontmatter}

\title{Scalability and Fragility in Bounded-Degree Consensus Networks } 

\thanks{This work was supported by the Swedish Research Council through grant 2016-00861 and by the Australian Research Council through grant DP190102859.}

\author[First]{Emma Tegling} 
\author[Second]{Richard H. Middleton} 
\author[Second]{Maria M. Seron}

\address[First]{Institute for Data, Systems and Society, Massachusetts Institute of Technology, Cambridge, MA 02139, USA {\tt (tegling@mit.edu)}.}
\address[Second]{School of Electrical Engineering and Computing, The University of Newcastle, Callaghan, NSW 2308, Australia {\tt (richard.middleton, maria.seron@newcastle.edu.au)}.}

\begin{abstract}                
%
We investigate the performance of linear consensus algorithms subject to a scaling of the underlying network size. Specifically, we model networked systems with \nth order integrator dynamics over families of undirected, weighted graphs with bounded nodal degrees. In such networks, the algebraic connectivity affects convergence rates, sensitivity, and, for high-order consensus ($n\ge 3$), stability properties. This connectivity scales unfavorably in network size, except in expander families, where consensus performs well regardless of network size. We show, however, that consensus over expander families is \emph{fragile} to a grounding of the network (resulting in leader-follower consensus). We show that grounding may deteriorate system performance by orders of magnitude in large networks, or cause instability in high-order consensus. Our results, which we illustrate through simulations, also point to a fundamental limitation to the scalability of consensus networks with leaders, which does not apply to leaderless networks. \vspace{-2mm}
\end{abstract}
\begin{keyword}
Distributed control; Large-scale systems; Robustness
\end{keyword}
\end{frontmatter}

\section{Introduction}
\label{sec:intro}
\vspace{-2mm}
Starting with early works on consensus problems over networks, there has been a great deal of interest in dynamic systems properties of classes of networked dynamic systems. Fundamental questions have been posed related to convergence~\citep{OlfatiSaber2004}, controllability~\citep{Olshevsky2014, Pasqualetti2014}, and performance~\citep{Bamieh2012, SiamiMotee2014}. In these cases, there have been several situations where poor dynamic behaviors can be observed in large networks. 

A particular component of poor behavior in families of consensus networks can be described as \emph{scale fragility} (or just \emph{fragility}) wherein stability properties are lost for large-scale networks in the family. Here, \cite{Stuedli2017} provide an analysis of certain classes of cyclic networks, and \cite{Tegling2019} point to such scale fragilities in networks with high-order dynamics. 
Even if stability can be maintained for all elements of a family of networks, it is desirable that the behavior be \emph{scalable}, that is, that performance (such as time constants and sensitivity properties) be uniform with respect to network size within the family. See, e.g., discussions in \cite{Lestas2006, Bamieh2012} and~\cite{Tegling2019b}.  

%
%

In both cases, \emph{fragility} and \emph{scalability}, it turns out that the \emph{algebraic connectivity}, as defined by certain eigenvalue properties of the network graph Laplacian, plays a crucial role. 
For example, in first order consensus problems, the algebraic connectivity is directly related to the slowest mode in the exponential convergence to consensus \citep{OlfatiSaber2004}. Sensitivity and lack of network coherence can also be attributed to the algebraic connectivity approaching zero as the network size grows~\citep{SiamiMotee2014,Tegling2019b}. 
It is therefore of interest to consider families of networks where the algebraic connectivity may be bounded away from zero, independent of the network size. 
At the same time, many applications require communications overheads to be modest. It is therefore relevant to enforce a uniform (with respect to network size) upper bound on the nodal degrees. This is also a key assumption in our present work.


The two objectives -- bounded nodal degrees yet well-behaved algebraic connectivity -- are reconciled only in so-called \emph{expander families}. 
Expander families are typically characterized through combinatorial conditions ensuring that the networks are sufficiently interconnected~\citep{Alon1986}. As one of this paper's results, we complement those conditions with a tractable algebraic characterization of non-expander families (of undirected, weighted graphs).  

Not surprisingly, the fact that consensus algorithms perform well over expander networks has been observed in earlier consensus literature. 
In particular, \citeauthor{OlfatiSaber2005} has showcased the fast convergence properties of consensus in small-world networks (2005), and Ramanujan graphs (graphs that maximize the algebraic connectivity) (2007). 
\cite{Kar2008} show that Ramanujan graphs optimize the convergence speed of distributed inference problems, and \cite{Li2009} discuss quantized consensus over expanders.  
To the best of our knowledge, however, an issue that has \emph{not} been observed in the literature is the {fragility} of these results towards a \emph{grounding} of the network. This is the focus of the present paper.

Grounding a network implies that the state at one of the nodes is fixed, and made independent of neighboring nodes. The terminology originates from electrical networks; in the context of consensus one often speaks of \emph{leader-follower consensus} since the grounded node acts as a leader for the remaining network. Leader-follower consensus is natural in many contexts, like platooning after a lead vehicle~\citep{Seiler2004}, slack bus control in DC networks~\citep{Andreasson2017ACC}, and pinning control~\citep{Chen2007}. It may, however, also arise inadvertently, if a local controller ceases to function, if a node one-sidedly disconnects from its neighbors, or through a malicious attack. Either way, the good performance that was achievable in expander networks is inevitably lost. 

The dynamics of grounded networked systems are described by a grounded graph Laplacian. Therefore, performance aspects, which in standard consensus depend on the algebraic connectivity, now instead depend on the slowest mode of the grounded Laplacian (here termed \grev). While the algebraic connectivity can stay bounded away from zero in bounded-degree networks, the grounded eigenvalue is shown to always decrease in network size (in undirected graphs). This is a fundamental difference between the two types of consensus dynamics -- and one we wish to pinpoint here as an important fragility. 

The scalability and fragility properties we discuss apply to consensus algorithms of various orders. In Section~\ref{sec:setup}, we therefore set up our problem as an \nth order consensus problem over (families of) undirected, weighted graphs. In Section~\ref{sec:l2section} we review the role of the algebraic connectivity, discuss its unfavorable scaling in bounded-degree networks, and introduce expander families together with a discussion of random graphs. We show in Section~\ref{sec:fragile} that consensus over expander networks is fragile to network grounding, which we demonstrate through simulations in Section~\ref{sec:examples}. We conclude in Section~\ref{sec:conclusions}.

\vspace{-1mm}
\section{Problem setup}
\label{sec:setup}
\vspace{-2mm}
%
%

\subsection{Network model and definitions}
\label{sec:defs}
\vspace{-3mm}
Consider a network described by a graph $\mathcal{G} = \{\mathcal{V}, \mathcal{E}\}$ with $N = |\mathcal{V}|$ nodes. The set $\mathcal{E} \subset \mathcal{V} \times \mathcal{V}$ contains the edges, each of which has an associated weight $w_{ij}>0$. 
Denote by $\Nn_i$ the neighbor set of node $i$ in $\Gg$, and define the 
(vertex) \emph{boundary} of a set $X\subset \Vv$ as $\partial X = \{j \in \bar{X} ~|~ (i,j) \in \Ee,~i \in X \}$, where $\bar{X} = V\backslash X$.
Going forward, we will also consider sequences, or \emph{families} of graphs, denoted~$\{\Gg_N\}$, in which the network size $N$ is increasing. Throughout this paper, we assume all graphs to be undirected and connected. 

Denote by $L$ the weighted graph Laplacian of $\Gg$, whose elements $L_{ij} = \sum_{k = 1,k\neq i}^N w_{ik}$ if $i = j$ and $L_{ij} = -w_{ij}$ otherwise. The Laplacian eigenvalues are denoted $\lambda_i$ (or~$\lambda_i(\Gg)$ where explicitness is needed) for $ i = 1,\ldots,N$ and are ordered so that $0 = \lambda_1 < \lambda_2 \le \ldots \le \lambda_N$. The smallest non-zero eigenvalue $\lambda_2(\Gg)$, also known as the Fiedler eigenvalue, is called the \emph{algebraic connectivity} of the graph~$\Gg$.

Each node $i \in \Vv$ has an associated state $x_i(t)  \in \mathbb{R}$, which is assumed to describe deviations from a desired setpoint. Its time derivatives are denoted according to $x_i^{(0)}(t) = x_i(t)$, $x_i^{(1)}(t) = \ddt x_i(t)  = \dot{x}_i(t)$, $x_i^{(2)}(t) = \frac{\mathrm{d}^2}{\mathrm{d}t^2} x_i(t)  = \ddot{x}_i(t)$ etc.

\vspace{-2mm}
\subsection{The \nth order consensus problem}
\label{sec:consensus}
\vspace{-3mm}
We will consider examples of networked systems with various order of the local dynamics. 
For generality, we therefore model the local dynamics at each node $i\in \Vv$ as an \nth order integrator:
\begin{align} \nonumber
\ddt x_i^{(0)}(t)  &= x_i^{(1)}(t)\\ \nonumber
& ~\vdots \\ \nonumber
\ddt x_i^{(n-2)}(t)  &= x_i^{(n-1)}(t) \\  \nonumber
\ddt x_i^{(n-1)}(t)  &=  u_i(t).
\end{align}
We consider the \nth order consensus algorithm:
\begin{equation}
\label{eq:consensuscompact}
u_i(t) = - \sum_{k = 0}^{n-1} a_k\sum_{j \in \mathcal{N}_i} w_{ij}(x_i^{(k)}(t) - x_j^{(k)}(t)),
\end{equation}
where $a_k$ are nonnegative fixed gains. Going forward, we will drop the time dependence in the notation. 

Defining the full state vector $\xi = [x^{(0)},x^{(1)},\ldots,x^{(n-1)}]^T$, we can write the system's closed-loop dynamics as
\begin{equation}
\label{eq:closedloop}
\ddt \xi = \underbrace{\begin{bmatrix}
0 & I_N & 0 & \cdots & 0\\
0 & 0 & I_N & \cdots &\vdots \\
0 & 0 & 0 & \ddots &\vdots \\
0 & 0 & 0 & \cdots &I_N \\
-a_0L & -a_1L & -a_2L & \cdots & -a_{\mathrm{n-1}}L
\end{bmatrix}}_{A} \xi,
\end{equation}
where the graph Laplacian~$L$ was defined in Section~\ref{sec:defs} and $I_N$ denotes the $N\times N$ identity matrix. This model adheres to the one considered in, e.g.,~\cite{Ren2007} and is a straightforward extension of the better known first- and second-order algorithms. 

\vspace{-2mm}
\subsection{Leader-follower consensus in grounded networks}
\vspace{-3mm}
{Grounding} the network, by fixing the state at one of the nodes, results in a {leader-follower consensus} algorithm. Provided the consensus algorithm converges, it does so to the state at the grounded node (the leader). 

Without loss of generality, assume that node 1 is grounded and let its state be $x_1 = \dot{x}_1 = \ldots=x_1^{n} \equiv 0$. The closed-loop dynamics for the remaining nodes can be written as
\begin{equation}
\label{eq:closedloopred}
\ddt \bar{\xi}= \underbrace{\begin{bmatrix}
0 & I_{N-1} & 0 & \cdots & 0\\
0 & 0 & I_{N-1} & \cdots &\vdots \\
0 & 0 & 0 & \ddots &\vdots \\
0 & 0 & 0 & \cdots & I_{N-1} \\
-a_0\bar{L} & -a_1\bar{L}  & -a_2\bar{L}  & \cdots & -a_{n-1}\bar{L} \\
\end{bmatrix}}_{\bar{\mathcal{A}}} \bar{\xi} ,
\end{equation}
where $\bar{L}$ is the \emph{grounded} Laplacian obtained by deleting the first row and column of $L$ and~$\bar{\xi}$ is obtained by removing the states at node~1. The eigenvalues of $\bar{L}$ are denoted $\bar{\lambda_i}$ (or $\bar{\lambda}_i(\Gg)$) and are numbered as $0 < \bar{\lambda}_1\le \ldots < \bar{\lambda}_{N-1}$.
We will be particularly interested in the smallest eigenvalue, $\bar{\lambda}_1$, which we will refer to as the \emph{\grev}.
%


\vspace{-2mm}
\subsection{Underlying assumptions}
\vspace{-3mm}
In the upcoming sections, we will discuss properties of the system~\eqref{eq:closedloop} pertaining to its performance and robustness subject to a scaling of the network size. The following underlying assumptions on the system will be important for that discussion.
\begin{ass}[Bounded neighborhoods]
\label{ass:q}
Each local controller can receive measurements from at most $q$ neighbors, where the number $q$ is fixed and independent of $N$. That is, $ |\mathcal{N}_i| \le q$ for all $i \in \mathcal{V}$. 
\end{ass}
\begin{ass}[Bounded edge weights]
\label{ass:gainsweights}
The graph's edge weights are bounded, i.e., $0<w_{\min} \le w_{ij}\le w_{\max}<\infty$ for all $(i,j) \in \mathcal{E}$. These bounds hold for every graph in a family~$\{\Gg_N\}$.
\end{ass}
\begin{ass}[Fixed and bounded gains]
\label{ass:gainsweights}
The system's gains are bounded, i.e., $a_k\le a_{\max} <\infty $ for all $k = 0,1,\ldots, n$. They are also fixed, meaning that they do not change if the underlying network graph changes. In particular, they are independent of $N$.
\end{ass}
\vspace{-2mm}
Together, Assumptions 1 and 2 imply that the graph's nodal degrees remain bounded, even if the number of nodes increases. Assumption 3 implies that the local controller tunings are not affected by such an increase.
\vspace{-1mm}
\section{Connectivity scaling and expanders}
\label{sec:l2section}
\vspace{-2mm}
The connectivity of the network graph plays an important role for the performance of the consensus algorithm. Here, we will focus on the \emph{algebraic connectivity}, quantified through the smallest non-zero eigenvalue, $\lambda_2$, of the graph Laplacian. In this section, we review some -- both well and lesser known -- results on its role in consensus problems of different orders. We also discuss the scaling of $\lambda_2$ as networks grow and focus on expander families, which have particularly good connectivity properties. 

\vspace{-3mm}
\subsection{The role of $\lambda_2$}
\label{sec:l2properties}
\vspace{-3mm}
\subsubsection{Convergence rate}
Consider a first order consensus algorithm ($n = 1$). The rate of convergence is determined by the algebraic connectivity according to
\begin{equation}
 ||x(t) - x^{\mathrm{avg}}||\le ||x(0)-x^{\mathrm{avg}}|| e^{-a_0 \lambda_2 t} ,
\label{eq:convrate}
\end{equation}\citep{OlfatiSaber2004}, where $x^{\mathrm{avg}} \!= \!(\sum_{i = 1}^Nx_i)/N$ (an invariant quantity).
This implies that the speed at which a state of consensus is reached is inversely related to the size of $\lambda_2$. 
\vspace{-2mm}
\subsubsection{Sensitivity}
Assume that the system in~\eqref{eq:closedloop} is subject to a disturbance input:
\[ \dot{\xi} = A \xi + d,\]
where $d\in \mathbb{R}^N$, and let us consider the deviation from the consensus subspace $y_i = x_i - x^{\mathrm{avg}}$ 
as a measure of the algorithm's performance.
Now, denote by $G$ the input-output system from disturbance $d$ to the performance output $y$. Then, for first-order consensus ($n = 1$) it holds
\begin{equation}
\label{eq:hinfnorm}
 ||G||_\infty = \frac{1}{a_0\lambda_2}
\end{equation}
\citep{SiamiMotee2014}. 
The \hinf norm of a system has several interpretations. The interpretation as an induced norm (or $\Ll_2$ gain) $||G||_\infty = \sup_{d\neq 0} ||y||_2/||d||_2$ is particularly useful to characterize sensitivity. The relation~\eqref{eq:hinfnorm} thus implies that the system may amplify certain disturbance signals by a gain that is inversely proportional to~$\lambda_2$.


For second-order consensus ($n = 2$), the expression in~\eqref{eq:hinfnorm} instead gives a tight lower bound on the \hinf norm (follows from Theorem 2 in \cite{Pirani2017}).


\vspace{-2mm}
\subsubsection{Stability}
Now, consider the consensus algorithm~\eqref{eq:consensuscompact} with $n\ge 3$. A necessary (but not sufficient) condition for system stability, i.e., convergence to consensus, is  
\begin{equation}
\label{eq:highordercond}
 \lambda_2 > \frac{a_{n-3}}{a_{n-1}a_{n-2}}
\end{equation}
\citep[Theorem 3]{Tegling2019}. 
The condition implies that if the algebraic connectivity $\lambda_2(\Gg_N)\rightarrow 0$ as ${N\rightarrow \infty}$ in a graph family, then stability cannot be upheld beyond a certain network size (note, the $a_k$ are fixed by Assumption~\ref{ass:gainsweights}). 
High order ($n\ge 3$) consensus therefore has a {scale fragility} in such families of graphs. 

\vspace{-2mm}
\subsubsection{Grounded networks}
If the network is grounded, the above properties depend on the \grev~$\bll$ 
instead of on $\lambda_2$. 
That is, for leader-follower consensus with $n = 1$, the convergence rate is given by~$\bll$ and the \hinf norm from a disturbance to control error is $||\bar{G}||_\infty = 1/(a_0\bll)$ ~\citep{Pirani2017}. For leader-follower consensus with $n\ge 3$ a necessary stability condition reads  $\bll > {a_{n-3}}/{(a_{n-1}a_{n-2})}$~\citep[Theorem 5]{Tegling2019}.


\vspace{-1mm}
\subsection{Scaling of connectivity}
\label{sec:scaling}
\vspace{-4mm}
The algebraic connectivity tends not to scale well with network size in bounded-degree networks, leading to a \emph{lack of scalability} of the consensus algorithm. More precisely, for families of graphs $\{\Gg_N\}$ that do not satisfy certain \emph{expansion} properties it holds $\lambda_2(\Gg_N) \rightarrow 0$ as $N\rightarrow \infty$. Before discussing those properties, we will provide a more tractable algebraic description of graph families in which indeed $\lambda_2(\Gg_N) \rightarrow 0$. 

\begin{figure}[t]
\centering
\includegraphics[width = 0.4\textwidth]{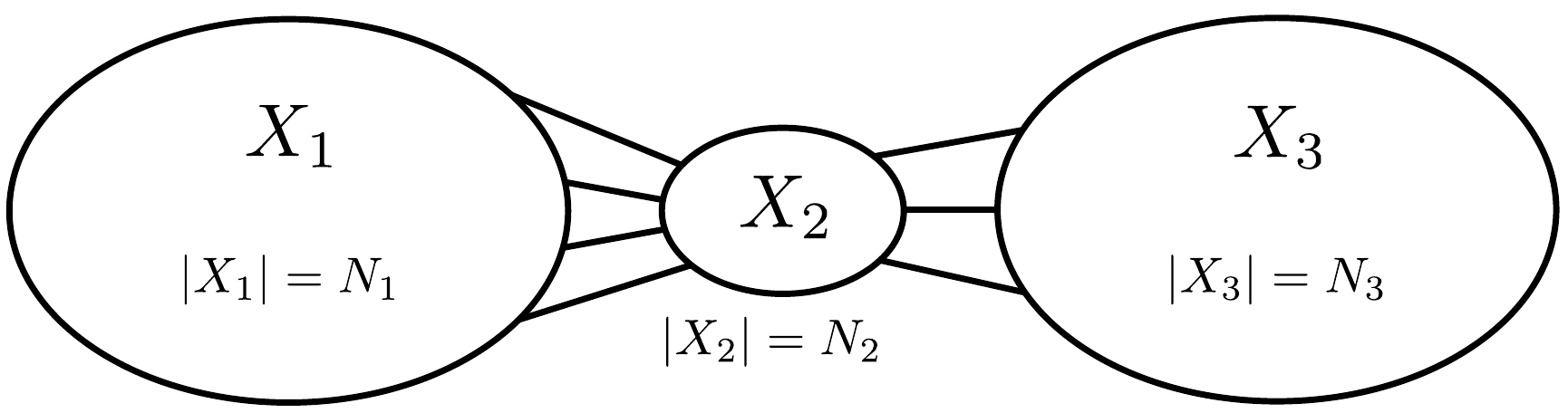}
\vspace{-2mm}
\caption{Partitioning of graph for Lemma~\ref{lem:bottleneck}. The set $X_2$ is a \emph{bottleneck} if it stays small compared to both $X_1$ and $X_3$ as the network grows. In this case, the algebraic connectivity decreases towards zero. } 

\label{fig:partition}
\end{figure}

For this purpose, partition a graph's vertex set into three disjoint sets $X_1,X_2,X_3$ so that $X_1\cup X_2 \cup X_3 = \Vv$ and ${|X_1| = N_1}$, ${|X_2| = N_2}$, ${|X_3| = N_3}$ as illustrated in Fig.~1. 
Each node in $X_2$ is connected to at least one node in both $X_1$ and $X_3$, but no edges connect $X_1$ and~$X_3$ directly.   
In other words, $X_2$ is the boundary set of both~$X_1$ and~$X_3$. This partitioning is always possible, unless the graph is complete (note that $X_1,X_2,X_3$ need not be connected subgraphs).  

By re-numbering the nodes, the graph Laplacian becomes
\begin{equation}
\Ll = \begin{bmatrix}
\Ll_1 & \Ll_{12} & 0_{N_1\times N_3} \\
\Ll_{12}^T & \Ll_2 & \Ll_{32}^T \\
0_{N_3\times N_1} &  \Ll_{32} & \Ll_3
\end{bmatrix}.
\label{eq:blocks}
\end{equation}
 If $N_2$ can be made small in relation to both $N_1$ and $N_3$, we say that the graph has a \emph{bottleneck}. The following lemma shows that if the bottleneck remains as the network grows, then $\lambda_2(\Gg_N) \rightarrow 0$.
 

\begin{lemma}
\label{lem:bottleneck}
Consider a graph family $\{\Gg_N\}$ and let Assumptions~1--2 hold. If every graph $\Gg_N$ in the family can be partitioned as outlined above in such a way that $N_2/N_1 \rightarrow 0$ and $N_2/N_3 \rightarrow 0$ as $N\rightarrow \infty$, then $\lambda_2(\Gg_N) \rightarrow 0$ as $N\rightarrow \infty$.
\end{lemma}
\vspace{-1mm}
\textit{Proof:} 
By the Rayleigh-Ritz theorem
\begin{equation}
\lambda_2 \le \frac{v^T\Ll v}{v^Tv},~~~\forall v \bot \mathbf{1},~v\neq 0,
\label{eq:RR}
\end{equation} 
since $\mathbf{1}$ is the eigenvector corresponding to $\lambda_1 =0$.
Let $\Ll$ be partitioned as in~\eqref{eq:blocks} and choose 
\[ v = \begin{bmatrix}
N_3\mathbf{1}_{N_1} \\ 0_{N_2} \\ -N_1 \mathbf{1}_{N_3}
\end{bmatrix}, \]
for which we verify $v^T\mathbf{1}_{N} = N_3N_1 - N_1N_3 = 0$ and $v^Tv = N_3^2N_1 + N_1^2N_3$.

Now,  
\(v^T\Ll v = N_3\mathbf{1}_{N_1}^T N_3 \Ll_{1}\mathbf{1}_{N_1} + 0 + N_1\mathbf{1}_{N_3}^T N_1 \Ll_{3}\mathbf{1}_{N_3} = N_3^2 d_{12} + N_1^2 d_{32},\)
where $d_{12}$ ($d_{32}$) is the total weight of all edges connecting $X_1$ ($X_3$) to $X_2$. By Assumptions~1--2 $d_{12},d_{32} \le q w_{\max}N_2$

Inserting in~\eqref{eq:RR} gives
\begin{equation}
\lambda_2 \le \frac{ N_3^2 d_{12} + N_1^2 d_{32}}{N_3^2N_1 + N_1^2N_3} \le q w_{\max}\frac{ N_3^2N_2  + N_1^2N_2 }{N_3^2N_1 + N_1^2N_3}.
\label{eq:l2bound}
\end{equation}
Assume without loss of generality that $N_2/N_1 \ge N_2/N_3$ for every $\Gg_N$. Then~\eqref{eq:l2bound} gives
\[\lambda_2 \le q w_{\max} \frac{N_2}{N_1}\cdot\frac{ N_3^2N_1  + N_1^2N_3 }{N_3^2N_1 + N_1^2N_3} = q w_{\max} \frac{N_2}{N_1}, \]
and, since $N_2/N_1 \rightarrow 0$ as $N\rightarrow \infty$, the result follows. \qed

Many families of bounded-degree graphs, for example planar graphs, lattices, and trees, will have bottlenecks in the sense of Lemma~\ref{lem:bottleneck}. They therefore have that $\lambda_2(\Gg_N) \rightarrow 0$ as $N\rightarrow \infty$. Such families were surveyed in~\cite{Tegling2019}, but Lemma~\ref{lem:bottleneck} provides a generalization. 

\begin{figure}
\centering
\includegraphics[width = 0.48\textwidth]{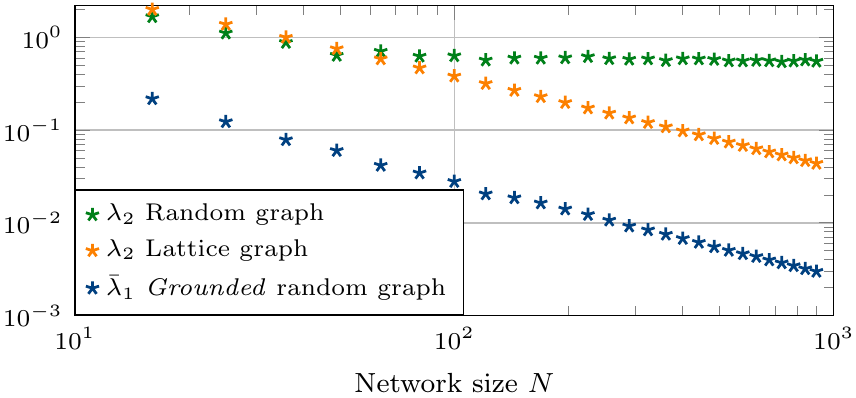}
\vspace{-3mm}
\caption{{Scaling of algebraic connectivity~$\lambda_2$ of randomly generated vs. 2D lattice graphs. Their nodal degrees are the same, yet $\lambda_2$ approaches zero in the lattice graphs but not in the random graphs. However, when the same random graphs are grounded, the \grev~$\bll$ approaches zero. }  }
\label{fig:l2scaling}
\end{figure}

 \vspace{-2mm}
\subsection{Expander families}
\vspace{-4mm}
The algebraic connectivity of a graph is related to its \emph{Cheeger constant}, defined here as 
\begin{equation}
h(\Gg) = \inf_{X \subset \Vv} \frac{ |\partial X|_d}{\min\{|X|_d, |\bar{X}|_d\} } ,
\label{eq:cheeger}
\end{equation}
where $\partial X$ is the boundary set of $X$ defined in Section~\ref{sec:defs}. Since we consider non-regular graphs, we measure a set~$W$ of nodes through $|W|_d := \sum_{i \in W}d_i$, where the nodal degree $d_i = \sum_{j \in \mathcal{N}_i} w_{ij}$ \citep{ChungBook}. 

Loosely speaking, a large Cheeger constant means that the graph has no bottlenecks; a partitioning as in Fig.~\ref{fig:partition} can only be done with a relatively large set~$X_2$. 
A precise definition of {expander families} follows:
\begin{definition}[Expander family]
Let $\{\Gg_N\}$ be a graph family in which $N\rightarrow \infty$. 
If the sequence $\{h(\Gg_N)\}$ is bounded away from zero, $\{ \Gg_N\}$ is an expander family. 
\end{definition}
\vspace{-2mm}
The following result is central to our discussion
\begin{result}
\label{res:expanders}
The graph family $\{ \Gg_N\}$ is an expander family if and only if the sequence $\{\lambda_2(\Gg_N)\}$ is bounded away from zero as $N\rightarrow \infty$. 
\end{result}
\vspace{-2mm}
(See e.g. \cite{ChungBook, KrebsBook} for a proof.)
Result~\ref{res:expanders} implies that the algebraic connectivity can indeed be prevented from decreasing towards zero in bounded-degree networks, making the consensus algorithm scalable. The formal requirement is that the Cheeger constant in~\eqref{eq:cheeger} does not decrease towards zero. 

\vspace{-2mm}
\subsubsection{Explicit constructions and random graphs}
It is known that there exist $k$-regular expander families for every ${k\ge 3}$~\citep{Pinsker1973}, but constructing them is non-trivial. Research into explicit construction rules has been ongoing since the 1970's. We will not review such rules here, but refer interested readers to, e.g.,~\cite{KrebsBook}.

Interestingly though, random graphs will \emph{almost surely }be expander families for sufficiently large $N$~\citep{Friedman1991}. Here, we mean graphs that are random in the sense that edges are selected through equally likely permutations of the node set. The good expansion properties of certain random graphs are  exploited in~\cite{OlfatiSaber2007b, Kar2008}, as well as in our upcoming simulations. 

For networked control problems, this implies that a way to achieve good connectivity, and thus scalability, properties, is to assign neighbors randomly across the network. Bandwidth can stay limited, since the number of neighbors remains bounded, but an ability to communicate across the scale of the network would be required. A more detailed study on the construction of expander families from a practical perspective is part of ongoing work. 

\begin{figure}[t]
\centering
\includegraphics[width = 0.36\textwidth]{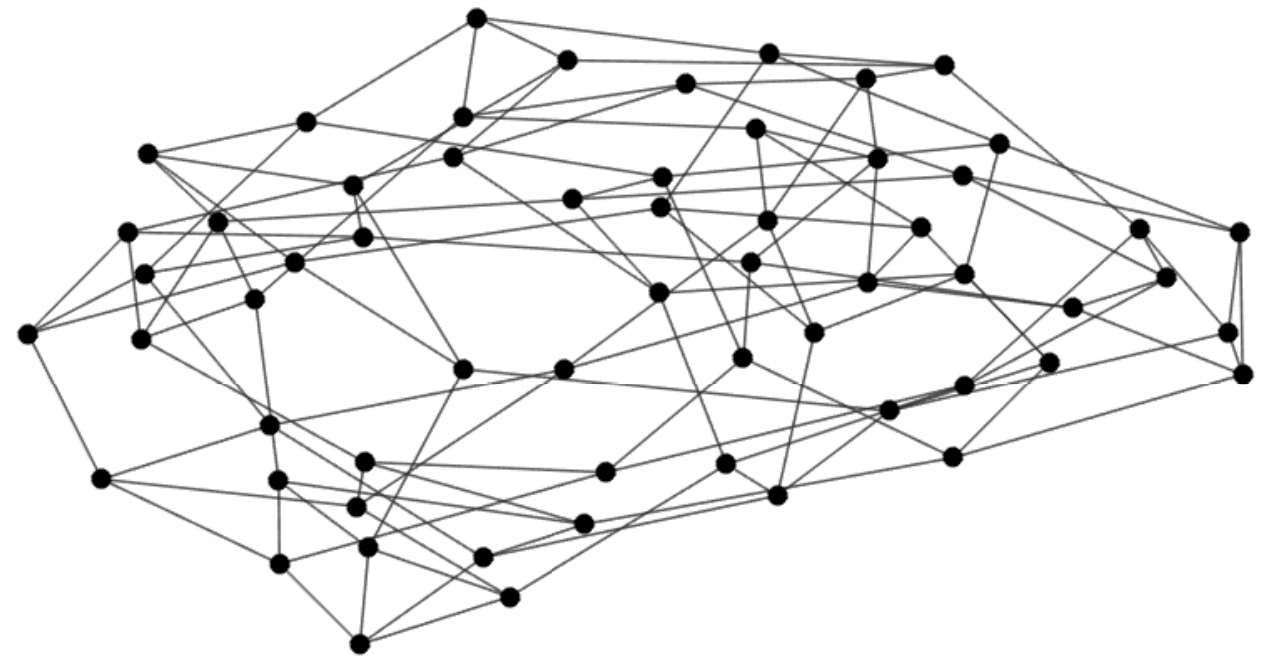}
\vspace{-2.5mm}
\caption{Example of a randomly generated graph as used in the examples. Here, $N = 60$ and each node has degree 4. With unit edge weights,  $\lambda_2\approx 0.64$ but $\bll  \approx 0.05$.} 
\label{fig:randomgraph}
\end{figure}

Fig.~\ref{fig:l2scaling} illustrates some of this section's results. Here, we have generated a sequence of random graphs of increasing size using the algorithm proposed by~\cite{Kim2006}. 
An example of one member of the sequence is given in Fig.~\ref{fig:randomgraph}.
Fig.~\ref{fig:l2scaling} compares their algebraic connectivity to that of 2-dimensional lattice graphs. 
In both cases, the degree of each node is~4, but their connectivities scale differently in~$N$. Fig.~\ref{fig:l2scaling} also displays the drastic discrepancy between $\lambda_2$ of the random graphs, and the corresponding \grev~$\bll$. 
This is the topic of the next section.



\section{Fragility towards network grounding}
\label{sec:fragile}
\vspace{-2mm}
In the previous section, we saw that maintaining good performance in growing consensus networks requires the underlying graphs to constitute an expander family. This can, for example, be achieved through random graphs. 
However, if the network is grounded, the \grev~$\bll$ \emph{inevitably} decreases towards zero as the network size grows. Consider the following lemma: 

\begin{lemma}
\label{lem:grounded}
Consider a graph family $\{\Gg_N\}$ and let Assumptions~1--2 hold. The smallest eigenvalue~$\bar{\lambda}_1(\Gg_N)$ of the grounded Laplacian~$\bar{\Ll}(\Gg_N)$ then satisfies
\begin{equation}
\label{eq:lambda1cond}
\bar{\lambda}_1(\Gg_N) \le \frac{q}{N-1}w_{\max}.
\end{equation}
\end{lemma}

\textit{Proof:}
By the Rayleigh-Ritz theorem it holds
\begin{equation*}
\bar{\lambda}_1 \le  \frac{v^T\bar{\Ll}v}{v^Tv} ,~~ \forall v\neq 0.
\end{equation*}
This implies in particular that 
\[ \bar{\lambda}_1 \le  \frac{\mathbf{1}_{N-1}^T\bar{\Ll}\mathbf{1}_{N-1}}{\mathbf{1}_{N-1}^T\mathbf{1}_{N-1}} = \frac{\sum_{k \in \mathcal{N}_1} w_{1k}}{N-1} \le \frac{q w_{\max}}{N-1},\]
where $\sum_{k \in \mathcal{N}_1} w_{1k}$ is the total weight of the edges leading to the grounded node 1. The equality holds since each row~$k$ of~$\bar{\Ll}$ sums to zero if the corresponding node~$k$ has no connection to the leader, and otherwise to~${w_{1k}\le w_{\max}}$. \qed

Lemma~\ref{lem:grounded} says that, under the given assumptions, $\bll \rightarrow 0$ as $N\rightarrow \infty$. Therefore, the performance of the leader-follower consensus algorithm never scales well in grounded bounded-degree networks. We next discuss some implications in more detail. 

\vspace{-2mm}
\subsection{Implications}
\vspace{-3mm}
\subsubsection{Performance degradation}
Lemma~\ref{lem:grounded} implies that the consensus algorithm can be \emph{fragile} to a grounding of the network. Consider a scenario where a large bounded-degree network has been carefully designed to avoid the bottlenecks from Lemma~\ref{lem:bottleneck} to ensure $\lambda_2$ is large. For example, the network in Fig.~\ref{fig:randomgraph}. Assume first that the consensus algorithm~\eqref{eq:consensuscompact} with $n = 1$ is run over this network. If a single node (say, number 1) turns off its controller so that $u_1 = 0$, the system instead obeys the leader-follower dynamics~\eqref{eq:closedloopred}. Since we may have ${\bll << \lambda_2}$, the convergence time and sensitivity can increase radically. 

High-order ($n\ge 3$) consensus is yet more fragile. In this case, the breakdown of one controller, or the active decision of one agent to disconnect from its neighbors, would cause a grounding of the network. If the network is sufficiently large, the fact that $\bll << \lambda_2$ leads to a loss of stability. This scenario is simulated in Fig.~\ref{fig:third}.

\vspace{-2mm}
\subsubsection{Lack of scalability} Lemma~\ref{lem:grounded} and Result~\ref{res:expanders} together imply that there is an important difference between standard consensus and leader-follower consensus algorithms in their scalability properties. It is possible to achieve good scalability (in terms of the properties discussed in Section~\ref{sec:l2properties}) in standard consensus over bounded-degree networks, but it is fundamentally impossible in leader-follower consensus. 

The scalability of consensus in bounded-degree networks is therefore, in a sense, fragile to the assumption that the network has no leader. This has implications for, e.g., vehicle platooning problems. Here, one may wish to add communication links in an optimal way to increase connectivity and thereby improve performance (see, e.g.,~\cite{Swaroop2019}). When each vehicle has a bounded number of links, however, this cannot give a fully scalable performance as long as the platoon has an independent lead vehicle.

\begin{figure}
\includegraphics[width = 0.5\textwidth]{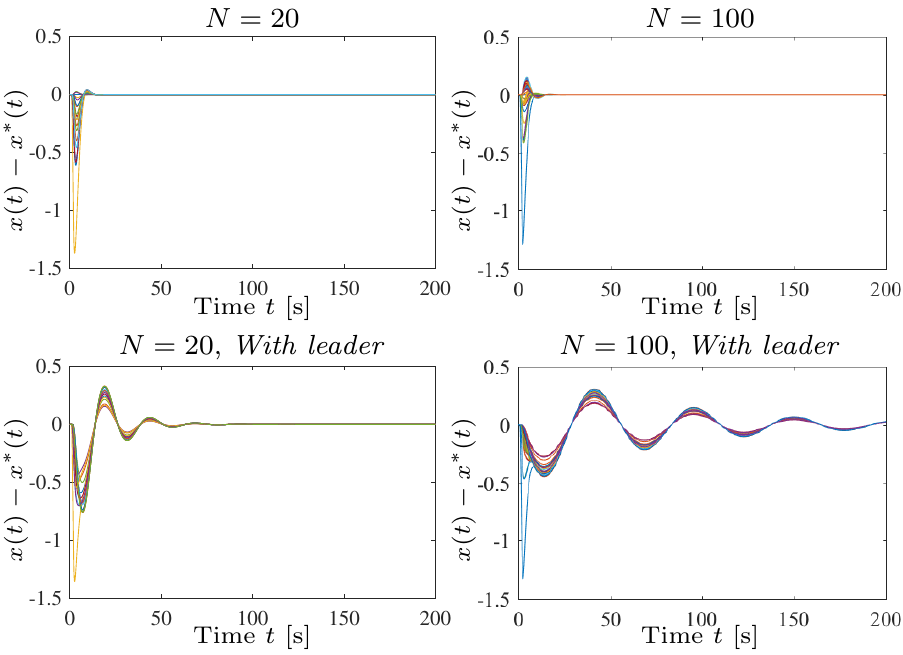}
\vspace{-3mm}
	\caption{Simulation of vehicular formations with and without leader, perturbed by a sudden deceleration of one vehicle. The system returns to the desired trajectory $x^*$ much faster in the absence of a lead vehicle (top panels), and the performance is similar in the smaller (left) and the larger  system (right). With a lead vehicle (bottom panels), scalability is limited.   }
	\label{fig:second}
\end{figure}

\vspace{-0.3mm}
\section{Examples}
\label{sec:examples}
\vspace{-3mm}
\subsection{Scalability of formation control }
\vspace{-4mm}
A cooperative vehicular formation control problem can be described through the consensus dynamics~\eqref{eq:closedloop} with $n = 2$. Here the objective is to steer $N$ vehicles to a trajectory where $x^*_i(t) = {v}^*t+\delta_i(t)$, with $v^*$ being a desired common velocity and $\delta_i$ a setpoint for each vehicle $i = 1,\ldots,N$. For example, in vehicle platooning, $v^*$ is a cruising velocity and $\delta_i$ is chosen to create desired spacings between vehicles.  

To illustrate the possible significant difference in performance scaling of the formation control problem with and without a lead vehicle, we simulate the system over randomly generated graphs with, respectively, $N=20$ and $N = 100$ nodes. Again, we use the algorithm proposed by~\cite{Kim2006} to generate the graphs. 
At time $t = 1$ s, one of the vehicles decelerates. As shown in Fig.~\ref{fig:second}, this disturbance is attenuated by the network in less than 20 seconds for both the smaller and larger leaderless networks. When the system has an independent lead vehicle, however, the network becomes grounded. It then takes near 4 times longer for the smaller network and about 20 times longer for the larger one (the full transient is not displayed). 
%
%




\vspace{-2mm}
\subsection{Loss of stability in third order consensus}
\begin{figure}
	\centering
	\includegraphics[scale=1]{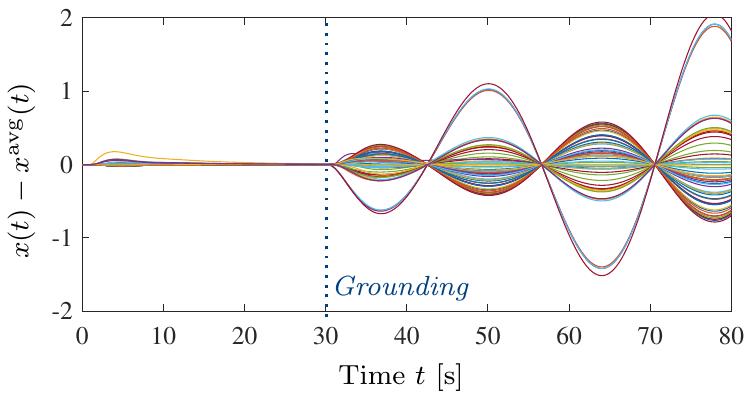}
	\vspace{-3mm}
	\caption{Position trajectories (relative to average) in simulation of $3^{\text{rd}}$ order consensus over the network in Fig.~\ref{fig:randomgraph}. After the network has been grounded at $t = 30$~s, the system is no longer stable and diverges when subjected to a disturbance -- a fragility. }
	\label{fig:third}
\end{figure}

\vspace{-4mm}
Consider the consensus algorithm~\eqref{eq:consensuscompact} with $n = 3$. We simulate the algorithm with $a_0 = 0.1$, $a_1 = 1$, $a_2 = 1$ over the graph displayed in Fig.~\ref{fig:randomgraph}, for simplicity with unit edge weights. At time $t = 1$ s, the network is subjected to a small acceleration disturbance at node~10. The disturbance is quickly attenuated. At time $t = 30$ s, the state of node~1 is fixed, that is, the network becomes grounded. The same disturbance is then applied at time $t = 31$ s. Now, the states diverge -- the system has become unstable. The position trajectories are displayed in Fig.~\ref{fig:third}.



\section{Conclusions}
\label{sec:conclusions}
\vspace{-2mm}
In this paper, we have considered classes of symmetric, high order consensus problems.
Looking at networks with uniformly bounded node degrees, we examined issues related to scalability. Specifically, the possibility of maintaining algebraic connectivity bounded away from zero as networks grow in size.
We have seen that the presence of `bottlenecks' in the network, or grounding (which creates leader-follower networks) guarantee that the algebraic connectivity (or the \grev, which plays a similar role) approaches zero as the network size increases. 

There are bounded-degree networks, `expander families', where algebraic connectivity scales well. Such families can (almost surely) be constructed through random graphs. However, their performance is highly fragile to network grounding, particularly in large networks, because the grounded eigenvalue is much smaller than the algebraic connectivity. 

Ongoing and future work include studies of, on one hand, design rules to achieve scalability in practical networked systems. On the other, strategies to mitigate the identified fragilities towards grounding. 

%

\bibliography{consrefs}   

\end{document}